\documentclass[12pt,a4paper]{article}
\usepackage[T1]{fontenc}
\usepackage[utf8]{inputenc}
\usepackage[english]{babel}
\usepackage{lmodern,microtype}
\usepackage[left=28mm,right=28mm,top=38mm,bottom=40mm,headheight=16pt,headsep=9mm]{geometry}
\usepackage{setspace}
\usepackage{amsmath,amssymb,amsthm,mathtools,mathrsfs}
\numberwithin{equation}{section}
\newtheoremstyle{plainspaced}{12pt}{12pt}{\itshape}{}{\bfseries}{.}{.6em}{}
\theoremstyle{plainspaced}
\usepackage{aliascnt}
\newtheorem{theorem}{Theorem}[section]
\newaliascnt{proposition}{theorem}
\newtheorem{proposition}[proposition]{Proposition}
\aliascntresetthe{proposition}
\newaliascnt{lemma}{theorem}
\newtheorem{lemma}[lemma]{Lemma}
\aliascntresetthe{lemma}
\newaliascnt{corollary}{theorem}

\aliascntresetthe{corollary}
\newtheoremstyle{uprightspaced}{12pt}{12pt}{\normalfont}{}{\bfseries}{.}{.6em}{}
\theoremstyle{uprightspaced}

\newtheorem{example}[theorem]{Example}

\theoremstyle{remark}
\newaliascnt{remark}{theorem}

\aliascntresetthe{remark}

\DeclareMathOperator{\id}{id}

\usepackage{titlesec}
\titleformat{\section}{\normalfont\large\scshape\centering}{\thesection.}{.65em}{}
\titlespacing*{\section}{0pt}{2.2\baselineskip}{.9\baselineskip}
\titleformat{\subsection}{\normalfont\normalsize\bfseries}{\thesubsection.}{.6em}{}
\titlespacing*{\subsection}{0pt}{1.4\baselineskip}{.5\baselineskip}
\usepackage{fancyhdr}
\newcommand{\shorttitle}{The Riemann--Roch formula away from the small primes}
\usepackage{enumitem}
\setlist[enumerate]{label=\textup{(\roman*)},leftmargin=*,topsep=.5em,itemsep=.3em}
\setlist[itemize]{leftmargin=*,topsep=.5em,itemsep=.3em}
\usepackage[hidelinks,hyperfootnotes=false]{hyperref}
\usepackage{titling}
\pretitle{\begin{center}\Large\bfseries}
\posttitle{\par\end{center}\vskip .7em}
\preauthor{\begin{center}\large\begin{tabular}[t]{c}}
\postauthor{\end{tabular}\par\end{center}}
\predate{\begin{center}\normalsize}
\postdate{\par\end{center}}

\usepackage{tikz-cd}
\usetikzlibrary{arrows.meta}
\newtheorem*{remark*}{Remark}
\allowdisplaybreaks[1]
\hypersetup{
  pdftitle={The Riemann--Roch formula away from the small primes},
  pdfsubject={Grothendieck--Riemann--Roch with localized integral coefficients},
  pdfkeywords={Grothendieck--Riemann--Roch, Chern character, Chow groups,
    algebraic K-theory, Adams operations, Todd classes},
  pdfdisplaydoctitle=true,
  bookmarksnumbered=true
}

\usepackage[nameinlink,noabbrev]{cleveref}
\crefname{theorem}{Theorem}{Theorems}
\crefname{proposition}{Proposition}{Propositions}
\crefname{lemma}{Lemma}{Lemmas}
\crefname{remark}{Remark}{Remarks}
\crefname{corollary}{Corollary}{Corollaries}
\crefname{section}{Section}{Sections}
\crefname{equation}{}{}
\creflabelformat{equation}{(#2#1#3)}

\newcommand{\CH}{\operatorname{CH}}
\newcommand{\Td}{\operatorname{Td}}
\newcommand{\ch}{\operatorname{ch}}
\newcommand{\Ch}{\operatorname{Ch}}
\newcommand{\cl}{\operatorname{cl}}
\newcommand{\gr}{\operatorname{gr}}
\newcommand{\rk}{\operatorname{rk}}
\newcommand{\Z}{\mathbf Z}
\newcommand{\Q}{\mathbf Q}
\newcommand{\Ftop}{F_{\mathrm{top}}}
\newcommand{\CHle}[1]{\CH^{\leq #1}}

\newcommand{\paperauthors}{A. Navarro Garmendia and J. Navarro Garmendia}
\newcommand{\paperaffiliations}{}
\newcommand{\paperemails}{%
\textsc{Alberto Navarro Garmendia}\\
Departamento de Matem\'atica Aplicada, Universidad Polit\'ecnica de Madrid,\\
Madrid, Spain\\
\href{mailto:alberto.navarro@upm.es}{\texttt{alberto.navarro@upm.es}}\\[1em]
\textsc{Jos\'e Navarro Garmendia}\\
Departamento de Matem\'aticas, Universidad de Extremadura,\\
Badajoz, Spain\\
\href{mailto:navarrogarmendia@unex.es}{\texttt{navarrogarmendia@unex.es}}}
\title{The Riemann--Roch formula \\ away from the small primes}
\author{\paperauthors}
\date{}
\hypersetup{pdfauthor={\paperauthors}}
\begin{document}
\begin{center}
{\Large\bfseries \MakeUppercase{\thetitle}\par}
\ifx\paperauthors\empty\else
  \vspace{1em}{\small\MakeUppercase{\paperauthors}\par}
\fi
\ifx\paperaffiliations\empty\else
  \vspace{.4em}{\small\paperaffiliations\par}
\fi

\end{center}
\thispagestyle{plain}

\begin{abstract}
Let  $f\colon X\to Y$ be a projective morphism between
smooth quasi-projective schemes of pure dimension over an arbitrary field, and put $d = \dim X - \dim Y$. For every $n\geq0$, we prove that the codimension-$n$ Grothendieck--Riemann--Roch formula holds in $\CH^n(Y)\otimes\Z[1/N!]$, where $N=\max\{n,n+d+1\}$. In particular, the formula retains the torsion information at every prime greater than $N$.
\end{abstract}

\begingroup
\renewcommand{\thefootnote}{}
\footnotetext{\textit{2020 Mathematics Subject Classification.}
Primary 14C40; Secondary 14C15, 19E20.\par
\textit{Key words and phrases.}
Grothendieck--Riemann--Roch, finite Chern character, Chow groups,
algebraic $K$-theory, Adams operations, Todd classes.}
\endgroup
\setcounter{footnote}{0}

\section{Introduction and main result}\label{sec:statement}

To motivate the problem we consider, let us first review  the classical Grothendieck--Riemann--Roch theorem~\cite{BorelSerre}. Let  $f\colon X\to Y$ be a projective morphism between
smooth quasi-projective schemes of pure dimension over a field $k$. 
We write
$f_!\colon K_0(X)\to K_0(Y)$ for the direct image in $K$-theory,
given by $f_![E]=\sum_j(-1)^j[R^jf_*E]$ under the identification
$K_0(Y)=G_0(Y)$, and $f_*$ for proper pushforward on Chow groups.

The Riemann--Roch formula, as generalized by Grothendieck, states that
\begin{equation}\label{eq:classical-grr}
 \ch(f_!\alpha)=f_*\bigl(\ch(\alpha)\Td(T_f)\bigr)
 \quad\text{in }\CH^*(Y)_{\Q} \ .
\end{equation}
Here
$T_f=T_X-f^*T_Y$ is the virtual relative tangent bundle, of rank
$d=\dim X-\dim Y$, and $\Td$ is the multiplicative characteristic
class associated with the series $z/(1-e^{-z})$.
The subscript indicates the coefficient ring: for an abelian group $A$
and a coefficient ring $B$, we write $A_B=A\otimes_{\Z}B$.
We use this convention for both Chow groups and $K$-groups.

Since $f_*$ lowers codimension by $d$, the degree-$n$ component of
\eqref{eq:classical-grr} is
\begin{equation}\label{eq:classical-grr-degree}
 \ch_n(f_!\alpha)=f_*\bigl[\ch(\alpha)\Td(T_f)\bigr]_{n+d}.
\end{equation}
We write $\ch_i=s_i/i!$, where $s_i$ is the $i$-th integral Newton
polynomial in the Chern classes for $i\geq1$, and $s_0=\rk$.
Brackets denote homogeneous components, with components of negative
degree understood to be zero.

Only finitely many denominators occur in
\eqref{eq:classical-grr-degree}, so that both sides can be defined
over a smaller subring of $\Q$. The question we address in this paper is whether they are
already equal over that ring, before passing to rational coefficients.

To specify the coefficient ring, recall that for an integer $r\geq0$
the Todd number
\[
 T_r=\prod_{p\ \mathrm{prime}}p^{\lfloor r/(p-1)\rfloor}
\]
is the least common denominator of the Todd polynomials through
degree $r$; see \cite[\S2]{Pappas} and the references there. A prime
divides $T_r$ if and only if $p\leq r+1$, so
\[
 R_r=\Z[1/T_r]=\Z[1/(r+1)!].
\]
Now fix $n\geq0$, put $m=n+d$,
and set $N=\max\{n,m+1\}$ and
\begin{equation}\label{eq:coefficient-rings}
 \Lambda=\Z[1/N!] \ ,
\end{equation}
so that $N=m+1$ and $\Lambda=R_m$ when $d\geq0$, while $N=n$ when $d<0$.
Both sides of \eqref{eq:classical-grr-degree} are defined over $\Lambda$.
When $m\geq0$, this localization is minimal for the universal polynomials
in independent Chern classes: the coefficient of $c_1^n$ in $\ch_n$ is
$1/n!$ for $n\geq1$, and the term multiplying $\rk(\alpha)$ on the
right is $\Td_m(T_f)$, whose universal denominator is $T_m$.
Together these denominators require exactly the primes $p\leq N$ to be
inverted, although for a particular morphism the formula may make sense
over a smaller ring.
Our main result, \cref{thm:single-square},
asserts that \eqref{eq:classical-grr-degree} holds over this ring, so that it retains torsion information at every prime $p>N$.

To express the result with a diagram (an image is worth a thousand words), we introduce
finite windows in $K$-theory and Chow groups. Suppose first that
$m=n+d\geq0$. Only degrees through $m$ on $X$ and through $n$ on $Y$
enter \eqref{eq:classical-grr-degree}. For $r\geq0$, define
\[
 K_{\leq r}(X)=K_0(X)/\Ftop^{r+1}K_0(X),
 \qquad \CHle r(X)_B=\bigoplus_{i=0}^r\CH^i(X)_B,
\]
where $\Ftop^{r+1}K_0(X)$ is generated by coherent-sheaf classes
supported in codimension at least $r+1$. These windows carry the
quotient topological filtration and the descending codimension
filtration, respectively, and the Chow product is truncated above
degree $r$. When $r!$ is invertible in a localization $B$ of $\Z$,
the finite Chern character is the map
\[
 \Ch_r=(s_i/i!)_{0\leq i\leq r}\colon
 K_{\leq r}(X)_B\longrightarrow\CHle r(X)_B.
\]
We show in \cref{prop:cycle-isomorphism} that this map is well defined
and is an isomorphism of filtered rings.

For a virtual bundle $E$ on $X$, we also write
\begin{equation}\label{eq:truncated-Todd}
 \Td_{\leq m}(E)=[\Td(E)]_{\leq m}\in\CHle m(X)_{R_m}.
\end{equation}
This is the multiplicative class defined by the Todd series through
degree $m$, taking quotients for differences of bundles.
The usual direct images $f_!$ and $f_*$ induce maps between the
corresponding windows, as recalled in \cref{subsec:direct-images}.
With this notation, our main result reads as follows. 

\begin{theorem}\label{thm:single-square}
Let $f\colon X\to Y$ be a projective morphism between smooth
quasi-projective schemes of pure dimension over a field, and put
$d=\dim X-\dim Y$. For $n\geq0$, set $m=n+d$,
$N=\max\{n,m+1\}$ and $\Lambda=\Z[1/N!]$.

If $m\geq0$, the following square commutes:
\begin{equation}\label{eq:single-square}
\begin{tikzcd}[column sep=7.2em,row sep=3em]
 K_{\leq m}(X)_{\Lambda}
   \arrow[r,"f_!"]\arrow[d,"\Ch_m"',"\sim"]
 &K_{\leq n}(Y)_{\Lambda}
   \arrow[d,"\Ch_n","\sim"']\\
 \CHle m(X)_{\Lambda}
   \arrow[r,"f_*(\Td_{\leq m}(T_f)\cdot(-))"']
 &\CHle n(Y)_{\Lambda}
\end{tikzcd}
\end{equation}
and its vertical arrows are isomorphisms of filtered rings.

In particular, taking the degree-$n$ component gives, for every
$\alpha\in K_0(X)$,
\begin{equation}\label{eq:main-relative}
 \ch_n(f_!\alpha)=f_*[\ch(\alpha)\Td(T_f)]_m
 \quad\text{in }\CH^n(Y)_{\Lambda}.
\end{equation}
The same formula holds when $m<0$, with both sides equal to zero.
\end{theorem}

\begin{example}\label{ex:mumford}
Let $f\colon X\to Y$ be a smooth projective morphism of smooth
quasi-projective schemes whose fibres are geometrically connected
curves of genus $g$. Let $E=f_*\omega_{X/Y}$ be the Hodge bundle and
put $\kappa_i=f_*\bigl(c_1(\omega_{X/Y})^{i+1}\bigr)$. Here $d=1$, so
$N=n+2$. Apply \cref{thm:single-square} to $\alpha=[\mathcal O_X]$.
By relative Serre duality, $R^1f_*\mathcal O_X\cong E^\vee$, so
$f_![\mathcal O_X]=1-[E^\vee]$, while
$T_f=\omega_{X/Y}^\vee$ is a line bundle, whose Todd series has no odd
terms beyond degree one. Writing $B_j$ for the Bernoulli numbers defined
by $z/(e^z-1)=\sum_{j\geq0}B_j z^j/j!$, we obtain from
\eqref{eq:main-relative}, for $k\geq1$,
\begin{align*}
 \ch_{2k-1}(E)&=\frac{B_{2k}}{(2k)!}\,\kappa_{2k-1}
   &&\text{in }\CH^{2k-1}(Y)_{\Z[1/(2k+1)!]},\\
 \ch_{2k}(E)&=0
   &&\text{in }\CH^{2k}(Y)_{\Z[1/(2k+2)!]}.
\end{align*}
These are Mumford's relations \cite{Mumford}, obtained there with
rational coefficients. For $k=1$, the first reads $12\,c_1(E)=\kappa_1$
in $\operatorname{Pic}(Y)_{\Z[1/6]}$. 

\end{example}

Thus, our formula falls within the tradition of the search for integral Riemann--Roch statements. This has a long history that starts with 
Grothendieck's Riemann--Roch theorem without denominators, proved by
Jouanolou for regular closed immersions \cite{Jouanolou} (see also
\cite[Theorem~15.3]{FultonIntersection}), and with landmarks such as the work of
Fulton--MacPherson on characteristic classes of direct images under
covering maps \cite{FultonMacPherson} or Pappas's integral
Grothendieck--Riemann--Roch theorem
\cite[Theorem~2.2 and Corollary~2.3]{Pappas}. 
The result most directly related to ours is Pappas's; in characteristic zero,
his formulas are stronger than \eqref{eq:main-relative}: they hold
in the integral Chow groups with explicit factors determined by the
Todd numbers, and imply our formula after inverting $T_M$, where
$M=\max\{m,n\}$. Pappas's proof uses resolution of singularities and
the weak factorization theorem for birational maps, and is consequently
restricted to characteristic zero. Our result gives a weaker form
of integrality, but one that holds over an arbitrary field.

In our companion paper \cite{NavarroNavarroHigher}, we give Riemann--Roch
formulas with explicit integer multipliers for direct images in higher
algebraic $K$-theory and Bloch's higher Chow groups. These formulas hold
integrally in characteristic zero and over $\Z[1/p]$ in characteristic
$p>0$. In another companion paper \cite[Theorem~1.1]{NavarroNavarroIntegral},
we use Adams operations to obtain integral Riemann--Roch formulas on $K_0$
over arbitrary fields, with explicit integer multipliers depending only
on the degrees involved.

As for the proof, the standard argument to derive \eqref{eq:classical-grr} factors $f$ as a
closed immersion $i\colon X\to\mathbf P^r_Y$ followed by a projection
$\pi\colon\mathbf P^r_Y\to Y$. Since $\pi_*$ lowers codimension by $r$,
the degree-$n$ formula for $\pi$ involves $\ch_{n+r}(i_!\alpha)$ on
$\mathbf P^r_Y$. There is no bound on $r$ in terms of $n$ and $d$, so
this intermediate Chern character may require primes that are not
inverted in $\Lambda$.
On the other hand, it is worth noting that the Riemann--Roch formula for the Adams operation $\psi^a$ is an
identity in $K$-theory in which only $a$ needs to be inverted, see \cref{eq:ARR}; although that theorem is also proved through factorization, the
auxiliary dimension $r$ leaves no trace on the coefficients of the formula. 

Our argument builds on an observation of Fulton--Lang: over $\Q$, the
Riemann--Roch formula for a single Adams operation $\psi^a$, with
$a\geq2$, determines the formula for the Chern character
\cite[Ch.~III, \S4]{FultonLang}; see also
\cite[Corollary~4.7]{NavarroNavarroAdams} for a diagrammatic description.
We apply this weight-separation argument to finite windows over
$\Lambda$. Since the finite
Chern character is an isomorphism over this ring (\cref{prop:cycle-isomorphism}), the main issue is arithmetic: over $\Q$, distinct weights can be separated by dividing by their
difference; over $\Lambda$, we need to use the operations $\psi^a$, with
$1\leq a\leq m+1$,w altogether (\cref{lem:unit-ideal} and
\cref{prop:weight-uniqueness}).

\section{The cycle map after a finite localization}\label{sec:finite-ch}

Let $X$ be a smooth quasi-projective scheme over $k$. We show that,
through degree $r$, the cycle map becomes an isomorphism after inverting
$r!$. Its inverse on the associated graded ring is induced by the finite
Chern character (\cref{prop:cycle-isomorphism}). We begin by recalling
our conventions for filtered maps and rings.

Our filtrations are finite, descending and multiplicative: for some
$\ell\geq0$,
\[
 A=F^0A\supseteq F^1A\supseteq\cdots\supseteq F^{\ell+1}A=0,
 \qquad F^iA\,F^jA\subseteq F^{i+j}A.
\]
We set $F^iA=A$ for $i\leq0$ and $F^iA=0$ for $i>\ell$.
An additive map $u\colon A\to A'$ is \emph{filtered of degree $e$}
if $u(F^iA)\subseteq F^{i+e}A'$ for every $i$.
For a graded ring $A=\bigoplus_{i\geq0}A^i$, we require
$1\in A^0$ and $A^iA^j\subseteq A^{i+j}$. An additive map between
graded rings is \emph{homogeneous of degree $e$} if it sends $A^i$
into $(A')^{i+e}$, with graded pieces of negative degree understood
to be zero.

The associated graded ring of a filtered ring is
$\gr_F A=\bigoplus_{i\geq0}F^iA/F^{i+1}A$.
A map filtered of degree $e$ induces a homogeneous map of degree $e$,
with components
\[
 \gr_e^i(u)\colon\gr_F^i A\longrightarrow\gr_F^{i+e}A'.
\]
Conversely, a graded ring with only finitely many nonzero pieces
carries the finite filtration $F^qA=\bigoplus_{i\geq q}A^i$.
Its associated graded ring is canonically identified with $A$.
Under this identification, a homogeneous map agrees with its
associated-graded map.

For $K_0(X)$, we use the topological filtration. It vanishes above
$\dim X$, its first step is the kernel of the
rank map,
and each step $\Ftop^qK_0(X)$ is an ideal, because tensoring with a
vector bundle does not enlarge supports. 

For $q\geq0$, the cycle map
\begin{equation}\label{eq:cycle-map}
 \cl_q^X\colon\CH^q(X)\longrightarrow\gr_{\Ftop}^qK_0(X)
       =\Ftop^qK_0(X)/\Ftop^{q+1}K_0(X)
\end{equation}
sends the class of an integral closed subscheme $Z\subset X$ of
codimension $q$ to $[\mathcal O_Z]$ modulo $\Ftop^{q+1}$.

The topological filtration is multiplicative, and the maps $\cl_q^X$
assemble to a surjective homomorphism of graded rings
\[
 \cl^X\colon \CH^*(X)\longrightarrow \gr_{\Ftop}K_0(X);
\]
see \cite[Appendix to Exp.~0, Ch.~II, Cor.~1 to Thm.~2.12]{SGA6}.

Nevertheless, the kernel can be nonzero;
see \cite[Example~5.14]{HautionMerkurjev}. The Riemann--Roch formula without denominators  (see \cite[Example~15.3.6]{FultonIntersection}) gives some information about this integral cycle
map. For
$\alpha\in\Ftop^qK_0(X)$ with $q\geq1$, it gives
\begin{equation}\label{eq:leading-chern}
 c_i(\alpha)=0\quad(0<i<q),\qquad
 c_q(\alpha)=(-1)^{q-1}(q-1)!b_q,
\end{equation}
where $b_q\in\CH^q(X)$ is any lift of
$[\alpha]\in\gr_{\Ftop}^qK_0(X)$ under $\cl_q$. 

Since $s_0(\alpha)=0$, Newton's identities then give
\begin{equation}\label{eq:leading-newton}
 s_i(\alpha)=0\quad(0\leq i<q),\qquad
 s_q(\alpha)=(-1)^{q-1}q\,c_q(\alpha)=q!b_q.
\end{equation}

For $q\geq1$, the kernel of $\cl_q^X$ is killed by $(q-1)!$.
Since $\cl_q^X$ is surjective, it becomes an isomorphism after
inverting $(q-1)!$. In particular, it is an integral isomorphism for
$q=1,2$; the same holds for $q=0$ by the rank map.

In what follows, we will be interested in finite windows, for both $K$-theory and Chow, with their corresponding filtrations: $K_{\leq r}(X)$ inherits
the quotient filtration, while $\CHle r(X)$ carries the codimension
filtration and the product truncated above degree $r$.

As for the cycle map, we write
$\cl_{\leq r}=\bigoplus_{q=0}^r\cl_q$ for the truncation through
degree $r$ and omit the
superscript when the scheme is clear.

\begin{proposition}\label{prop:cycle-isomorphism}
Let $r\geq0$, and let $B$ be a localization of $\Z$ in which $r!$
is invertible. The cycle map induces an isomorphism of graded rings
\begin{equation}\label{eq:finite-cycle}
 \cl_{\leq r}\otimes B\colon
 \CHle r(X)_B \xrightarrow{\ \sim\ }
 \gr_{\Ftop}\bigl(K_{\leq r}(X)_B\bigr)
\end{equation}
whose inverse in degree $q$ is induced by $s_q/q!$.

Moreover, the finite Chern character
\begin{equation}\label{eq:finite-ch}
 \Ch_r=(s_i/i!)_{0\leq i\leq r}\colon
 K_{\leq r}(X)_B \xrightarrow{\ \sim\ } \CHle r(X)_B
\end{equation}
is an isomorphism of filtered rings, and
$\gr(\Ch_r)=(\cl_{\leq r}\otimes B)^{-1}$.
\end{proposition}

\begin{proof}
Since $B$ is flat over $\Z$, the degree-$q$ part of
$\gr_{\Ftop}(K_{\leq r}(X)_B)$ is
$(\gr_{\Ftop}^q K_0(X))_B$ for $0\leq q\leq r$ and is zero
otherwise. Thus \eqref{eq:finite-cycle} is the truncation of the
integral cycle map after extension of scalars; in particular, it is
a surjective homomorphism of graded rings.

In degree zero, the componentwise rank map is inverse to $\cl_0$.
For $1\leq q\leq r$, equation~\eqref{eq:leading-newton} shows that
$s_q$ vanishes on $\Ftop^{q+1}K_0(X)$ and that the map it induces
on $\gr_{\Ftop}^q K_0(X)$ satisfies
\[
 \bar s_q\circ\cl_q=q!\,\id_{\CH^q(X)}.
\]
Since $q!$ is a unit in $B$ and $\cl_q$ is surjective, it follows
that $\cl_q\otimes B$ is an isomorphism with inverse induced by
$s_q/q!$, thus proving the first assertion.

On $K_0(X)_B$, the truncated Chern character is a ring homomorphism
to $\CHle r(X)_B$ with its truncated product
\cite[\S3.1, Lemma~6]{MaillotRoessler}. By
\eqref{eq:leading-newton}, its components of degree less than $q$
vanish on $(\Ftop^q K_0(X))_B$. Thus it preserves the filtrations
and vanishes on $(\Ftop^{r+1}K_0(X))_B$, so it descends to
\eqref{eq:finite-ch}. Its associated-graded map in degree $q$ is
induced by $s_q/q!$, which we have just identified with
$(\cl_q\otimes B)^{-1}$. Hence $\gr(\Ch_r)$ is an isomorphism.
Since both filtrations are finite, a standard argument using descending induction shows that $\Ch_r$ maps each filtration step isomorphically
onto the corresponding step.
\end{proof}

Over $\Q$, this statement recovers the classical Chern-character isomorphism; see \cite[Ch.~II, Corollary~8.9.1]{WeibelKbook}.

\section{The Riemann--Roch formula, revisited}
\label{subsec:direct-images}

In this section we recast the commutativity of the square
\eqref{eq:single-square} as an equality between two maps of Chow
windows. The first is the usual proper pushforward $f_*$. The second,
which we call
$P_f$, is obtained by
transporting $f_!$ to Chow groups along the finite Chern characters of
\cref{prop:cycle-isomorphism}, after dividing by the truncated Todd
class of $T_f$. Both maps are filtered of degree $-d$, and the main
point of the section is that they already agree on associated graded
groups (\cref{lem:graded-images}): this is just the compatibility of
the cycle map with proper pushforward. For a class of codimension $i$, the component
of $P_f$ in codimension $i-d$ is therefore $f_*$. Any discrepancy
can only occur in components of codimension $j>i-d$. Thus, \cref{thm:single-square} reduces to
showing that $P_f$ is homogeneous, a task that will be accomplished in  \cref{sec:proof,sec:Adams-RR}.

Until we return to the case $m<0$ at the end of the proof of
\cref{thm:single-square},
$f\colon X\to Y$ is a
projective morphism as in \cref{thm:single-square},
$d=\dim X-\dim Y$, and $n\geq0$ is chosen so that $m=n+d\geq0$.
We set $N=\max\{n,m+1\}$ and work over $\Lambda=\Z[1/N!]$. Note that
$m!$, $n!$ and $1,\ldots,m+1$ are units in $\Lambda$, and that
$R_m\subseteq\Lambda$. In particular, the finite Chern characters
$\Ch_m$ and $\Ch_n$ are defined over $\Lambda$, and so is the truncated
Todd class $\Td_{\leq m}(T_f)$.

Since the Chow pushforward is homogeneous of degree $-d$, it induces a
map
\[
 f_*\colon\CHle m(X)_{\Lambda}\longrightarrow\CHle n(Y)_{\Lambda}.
\]
Similarly, $f_!$ induces a map
$K_{\leq m}(X)_{\Lambda}\to K_{\leq n}(Y)_{\Lambda}$ that is
filtered of degree $-d$. Indeed, with the convention
$\Ftop^qK_0=K_0$ for $q\leq0$, one has
\begin{equation}\label{eq:support-shift}
 f_!(\Ftop^qK_0(X))\subseteq\Ftop^{q-d}K_0(Y),
\end{equation}
because the support of each higher direct image $R^jf_*\mathcal F$ is
contained in $f(\operatorname{Supp}\mathcal F)$, whose dimension is at
most that of $\operatorname{Supp}\mathcal F$.

We now recall the power-series calculus used to define the Todd class;
see \cite[Proposition~2.2, Remark~2.5 and Theorem~2.6]{NavarroNavarroAdams}.
Our first Chern classes are $c_1^K(L)=1-[L^\vee]$ and the usual
$c_1^{\CH}(L)$.

For $A=K_0(-)_B$ or $\CH^*(-)_B$, a series $h\in B[[z]]$ defines
an additive operation $h_+\colon K_0(X)\to A(X)$ by summing $h$
over the $A$-Chern roots of a bundle. If $h(0)$ is a unit,
multiplying instead defines the multiplicative characteristic class
$h_\times$,
with $h_\times(E-F)=h_\times(E)h_\times(F)^{-1}$.
These constructions are natural and determined by their values
on line bundles. The additive extension is a unital ring morphism if
\[
 h(0)=1,\qquad h(F_A(x,y))=h(x)h(y),
\]
where $F_K(x,y)=x+y-xy$ and $F_{\CH}(x,y)=x+y$.
The same calculus applies to truncated Chow rings, using identities
modulo terms of total degree greater than the truncation bound.

Universal identities between these truncated series can be checked
in $\Q[[z]]$: once their coefficients lie in a localization $B$
of $\Z$, the injection $B[z]/(z^{r+1})\hookrightarrow
\Q[z]/(z^{r+1})$ proves them over $B$, before evaluation in Chow.

The Todd series is the inverse of $(1-e^{-z})/z$. Through degree $r$,
the latter is
\[
 \left[\frac{1-e^{-z}}z\right]_{\leq r}
    =\sum_{j=0}^r\frac{(-1)^jz^j}{(j+1)!}\in R_r[z],
\]
a polynomial with constant term $1$; hence its inverse modulo $z^{r+1}$
also has coefficients in $R_r$. The multiplicative extension of this inverse defines the Todd class
$\Td_{\leq r}(E)\in\CHle r(X)_{R_r}$, for every virtual bundle $E$.

The class $\Td_{\leq m}(T_f)$ is a unit: it has the form $1+v$
with $v$ in the nilpotent positive-degree ideal.
By \cref{prop:cycle-isomorphism}, we may therefore define
\begin{equation}\label{eq:two-images}
 \begin{aligned}
 P_f&\colon\CHle m(X)_{\Lambda}\longrightarrow\CHle n(Y)_{\Lambda},\\
 P_f(u)&:=\Ch_n f_!\Ch_m^{-1}\bigl(\Td_{\leq m}(T_f)^{-1}u\bigr) , 
 \end{aligned}
\end{equation}
and reformulate the desired square as an equality of direct images:
\begin{lemma}\label{lem:normalized-square}
The square \eqref{eq:single-square} commutes if and only if $P_f=f_*$.
\end{lemma}
\begin{proof}
Every $u\in\CHle m(X)_{\Lambda}$ has the unique form
$u=\Td_{\leq m}(T_f)\Ch_m(\alpha)$ with
$\alpha\in K_{\leq m}(X)_{\Lambda}$. Evaluating $P_f=f_*$ at these
elements gives exactly \eqref{eq:single-square}.
\end{proof}

But we already know many things about this direct image $P_f$: it is clearly filtered of degree $-d$ (because multiplication by $\Td_{\leq m}(T_f)^{-1}$ preserves the source filtration, the finite
Chern characters are filtered isomorphisms, and $f_!$ satisfies
\eqref{eq:support-shift}) and its associated-graded map is
what it should be:

\begin{lemma}\label{lem:graded-images}
The associated-graded
map of $P_f$ is the usual Chow pushforward; that is,
\begin{equation}\label{eq:graded-images}
 \gr_{-d}^{\,i}(P_f)=f_*\colon
 \CH^i(X)_{\Lambda}\longrightarrow\CH^{i-d}(Y)_{\Lambda}
\end{equation}
for any $0\leq i\leq m$.
\end{lemma}

\begin{proof}
Multiplication by $\Td_{\leq m}(T_f)^{-1}$ induces the identity on
the associated graded ring, since its degree-zero component is $1$.
If $i-d<0$, both sides of \eqref{eq:graded-images} vanish, because
the target group is zero. Otherwise $0\leq i-d\leq n$, and
\cref{prop:cycle-isomorphism} gives
\begin{align*}
 \gr_{-d}^{\,i}(P_f)
 &=\gr_{-d}^{\,i}(\Ch_n f_!\Ch_m^{-1})\\
 &=(\cl_{i-d}^Y\otimes \Lambda)^{-1}
   \circ\gr_{-d}^{\,i}(f_!)\circ(\cl_i^X\otimes \Lambda).
\end{align*}
On the other hand, the cycle map commutes with proper pushforward
\cite[\S3.1]{HautionMerkurjev}; that is, for $u\in\CH^i(X)$,
\begin{equation}\label{eq:cycle-pushforward}
 \gr_{-d}^{\,i}(f_!)\bigl(\cl_i^X(u)\bigr)
       =\cl_{i-d}^Y(f_*u),
\end{equation}
and substituting this into the previous equality proves
\eqref{eq:graded-images}.
\end{proof}

The strategy now is clear: if we prove that $P_f$ is homogeneous, we are done, as it would then equal its
associated-graded map, $f_*$.

\section{Compatibility with the Adams operations}\label{sec:proof}

To prove that $P_f$ is homogeneous, the idea is to detect homogeneity through weights: on a Chow window,
let $\Phi^a$ act on the codimension-$i$ part as multiplication by
$a^i$, so that a homogeneous map of degree $-d$ intertwines these operators on source and target up to the scalar $a^{-d}$ (see \cref{eq:common-homogeneous}). Conversely, a map that is already filtered of degree $-d$ is homogeneous
if it satisfies this identity for every $a=1,\ldots,m+1$
(\cref{prop:weight-uniqueness}).

For a positive integer $a$, let $\Phi_X^a$ be the ring endomorphism
of the Chow window given by $\Phi_X^a|_{\CH^i(X)_{\Lambda}}=a^i\id$.
If $a$ is a unit in $\Lambda$, every map $u$ homogeneous of degree $-d$
(for example, $f_*$) satisfies
\begin{equation}\label{eq:common-homogeneous}
 \Phi_Y^a u=a^{-d}u\Phi_X^a.
\end{equation}
The key input is that $P_f$ satisfies the same identity:

\begin{proposition}\label{prop:common-Adams}
For every positive integer $a$ invertible in $\Lambda$, the normalized
direct image $P_f$ satisfies
\begin{equation}\label{eq:common-Adams}
 \Phi_Y^aP_f=a^{-d}P_f\Phi_X^a.
\end{equation}
\end{proposition}

In \cref{sec:Adams-RR} we prove this proposition by writing the
Adams--Riemann--Roch theorem in Chow coordinates. Here we use it to
prove homogeneity, starting with the arithmetic fact that allows us to
separate the weights over $\Lambda$.

\begin{lemma}\label{lem:unit-ideal}
Let $r\geq0$ and recall that $R_r=\Z[1/(r+1)!]$.
For distinct $i,j\in\{0,\ldots,r\}$, the differences
$a^j-a^i$, with $1\leq a\leq r+1$, generate the unit ideal:
\begin{equation}\label{eq:weight-unit-ideal}
 \bigl(a^j-a^i\mid 1\leq a\leq r+1\bigr)=R_r.
\end{equation}
Consequently, an element of an $R_r$-module annihilated by all
these differences is zero.
\end{lemma}

\begin{proof}
Suppose that the ideal in \eqref{eq:weight-unit-ideal} were proper.
It would then be contained in a maximal ideal $pR_r$, for some prime
$p>r+1$. The polynomial $t^j-t^i\in\mathbf F_p[t]$ is nonzero and has
degree at most $r$. However, the containment would make each of
$1,\ldots,r+1$ a root of it, and these are $r+1$ distinct elements of
$\mathbf F_p$, a contradiction. Thus the ideal is the whole ring, and
there exist $u_1,\ldots,u_{r+1}\in R_r$ with
\begin{equation}\label{eq:weight-Bezout}
 \sum_{a=1}^{r+1}u_a(a^j-a^i)=1.
\end{equation}
If all the differences annihilate an element $v$ of an $R_r$-module,
multiplying \eqref{eq:weight-Bezout} by $v$ gives $v=0$.
\end{proof}

\begin{proposition}\label{prop:weight-uniqueness}
Let $P\colon\CHle m(X)_{\Lambda}\longrightarrow\CHle n(Y)_{\Lambda}$
be $\Lambda$-linear and filtered of degree $-d$. If, moreover,
\begin{equation}\label{eq:uniqueness-Adams}
 \Phi_Y^a P=a^{-d}P\Phi_X^a\qquad(1\leq a\leq m+1),
\end{equation}
then $P$ is homogeneous of degree $-d$.
\end{proposition}

\begin{proof}
Let $u\in\CH^i(X)_{\Lambda}$, and write
$P(u)=\sum_{j=0}^n v_j$ with $v_j\in\CH^j(Y)_{\Lambda}$.
Since $P$ is filtered of degree $-d$, we have $v_j=0$ for
$j<i-d$. On the other hand, \eqref{eq:uniqueness-Adams} gives
\[
 \Phi_Y^aP(u)=a^{i-d}P(u)\qquad(1\leq a\leq m+1),
\]
so that $(a^j-a^{i-d})v_j=0$ for every $j$. Each such $a$ is a unit in
$\Lambda$, so for $j>i-d$ this becomes
\begin{equation}\label{eq:error-weight-gap}
 (a^\delta-1)v_j=0,
 \qquad\delta=j+d-i,
 \qquad1\leq\delta\leq n+d-i=m-i\leq m.
\end{equation}
(When $m=0$, this bound shows that there is no such $j$.)
By \cref{lem:unit-ideal}, applied with $r=m$ and weights $0$ and
$\delta$, these differences generate the unit ideal in $R_m$, and
hence in $\Lambda\supseteq R_m$. Therefore $v_j=0$ for $j>i-d$ as well.
Only the component of degree $i-d$ can survive, and if this degree
is negative, all components vanish. This proves homogeneity.
\end{proof}

\begin{proof}[Proof of \cref{thm:single-square}]
Suppose $m\geq0$. The map $P_f$ is filtered of degree $-d$, and
its associated-graded map is $f_*$ by \cref{lem:graded-images}.
Since $1,\ldots,m+1$ are units in $\Lambda$,
\cref{prop:common-Adams,prop:weight-uniqueness} show that $P_f$
is homogeneous of degree $-d$. It therefore equals its
associated-graded map, namely $f_*$. By \cref{lem:normalized-square},
this proves the commutativity of
\eqref{eq:single-square}. Its vertical arrows are filtered ring
isomorphisms by \cref{prop:cycle-isomorphism}.

Taking the degree-$n$ component gives \eqref{eq:main-relative},
since $f_*$ lowers degree by $d$ and $m=n+d$.

If $m<0$, \eqref{eq:support-shift} gives
$f_!\alpha\in\Ftop^{-d}K_0(Y)\subseteq\Ftop^{n+1}K_0(Y)$.
Hence $\ch_n(f_!\alpha)=0$ in $\CH^n(Y)_{\Lambda}$ by
\eqref{eq:leading-newton}, since $m+1\leq0$ gives $\Lambda=\Z[1/n!]$; the right-hand side of
\eqref{eq:main-relative} is zero by convention.
\end{proof}

\section{Adams--Riemann--Roch in Chow coordinates}\label{sec:Adams-RR}

This section is devoted to proving \cref{prop:common-Adams}. In essence, that proposition is the usual Adams--Riemann--Roch theorem, but translated into the Chow window via the Chern character.

For a positive integer $a$, the Adams operation
$\psi^a\colon K_0(X)\to K_0(X)$ is the additive extension of
$(1-z)^{-a}$. It is a ring endomorphism characterized
by $\psi^a[L]=[L^{\otimes a}]$; see
\cite[Example~2.7.1]{NavarroNavarroAdams}.
The integral identity $s_i\psi^a=a^i s_i$ follows by evaluation on
line bundles and additive extension. Consequently, whenever $r!$ is
invertible in the coefficient ring,
\begin{equation}\label{eq:ch-adams}
 \Ch_r\psi^a=\Phi_X^a\Ch_r.
\end{equation}
For $r\in\{m,n\}$, the kernel of $\Ch_r$ on $K_0(X)_{\Lambda}$ is
$(\Ftop^{r+1}K_0(X))_{\Lambda}$ by \cref{prop:cycle-isomorphism}, so
\eqref{eq:ch-adams} shows that $\psi^a$ preserves this kernel and
descends to the finite window $K_{\leq r}(X)_{\Lambda}$. The induced
operation preserves the filtration, again by
\cref{prop:cycle-isomorphism} and \eqref{eq:ch-adams}.

Bott's cannibalistic class $\theta^a$ is determined by multiplicativity and
\[
 \theta^a(L)=1+[L]+\cdots+[L]^{a-1}
\]
for line bundles $L$. After inverting $a$, it extends to virtual
bundles by $\theta^a(E-F)=\theta^a(E)\theta^a(F)^{-1}$.
Write $T_f^\vee=[\Omega^1_{X/k}]-f^*[\Omega^1_{Y/k}]\in K_0(X)$
for the virtual cotangent class.
The Adams--Riemann--Roch theorem states that
\begin{equation}\label{eq:ARR}
 \psi_Y^a f_!(\alpha)
 =f_!\bigl(\theta^a(T_f^\vee)^{-1}\psi_X^a(\alpha)\bigr)
 \quad\text{in }K_0(Y)_{\Z[1/a]};
\end{equation}
see \cite[Ch.~V, Theorem~7.6]{FultonLang} and
\cite[Corollary~4.5]{NavarroNavarroAdams}.
For $a$ invertible in $\Lambda$, this descends to the commutative square
\begin{equation}\label{eq:arr-square}
\begin{tikzcd}[column sep=large,row sep=large]
 K_{\leq m}(X)_{\Lambda}
   \arrow[r,"f_!"]\arrow[d,"\theta^a(T_f^\vee)^{-1}\psi_X^a"']
 &K_{\leq n}(Y)_{\Lambda}\arrow[d,"\psi_Y^a"]\\
 K_{\leq m}(X)_{\Lambda}\arrow[r,"f_!"']
 &K_{\leq n}(Y)_{\Lambda}.
\end{tikzcd}
\end{equation}
Indeed, the Adams operations preserve the window ideals by the previous
paragraph, multiplication by the Bott class preserves them because they
are ideals, and $f_!$ respects them by \eqref{eq:support-shift}.

Let us write $t=\Td_{\leq m}(T_f)$ and consider the filtered
$\Lambda$-module isomorphisms
\[
 \begin{aligned}
 \rho_X\colon\CHle m(X)_{\Lambda}&\longrightarrow K_{\leq m}(X)_{\Lambda},
 &u&\longmapsto\Ch_m^{-1}(t^{-1}u),\\
 \rho_Y\colon\CHle n(Y)_{\Lambda}&\longrightarrow K_{\leq n}(Y)_{\Lambda},
 &v&\longmapsto\Ch_n^{-1}(v).
 \end{aligned}
\]
Their inverses are $\rho_X^{-1}(\alpha)=t\Ch_m(\alpha)$ and
$\rho_Y^{-1}=\Ch_n$, and, by the definition \eqref{eq:two-images},
\begin{equation}\label{eq:def-P}
 P_f=\rho_Y^{-1}f_!\,\rho_X.
\end{equation}

\begin{proof}[Proof of \cref{prop:common-Adams}]
Write $A_f^a=\theta^a(T_f^\vee)^{-1}\psi_X^a$ for the left vertical
arrow of \eqref{eq:arr-square}. We first observe that
\[
 \rho_X^{-1}A_f^a\rho_X=a^{-d}\Phi_X^a,
 \qquad \rho_Y^{-1}\psi_Y^a\rho_Y=\Phi_Y^a.
\]
The first equality is \cref{lem:Adams-Todd}, proved below; the second
is \eqref{eq:ch-adams}. Together with \eqref{eq:def-P}, they give the
four side faces of the following cube. The upper face is the
Adams--Riemann--Roch square \eqref{eq:arr-square} in $K$-theory.
The lower face expresses the same identity in Chow coordinates: on the
target the change of coordinates is $\Ch_n$, and on the source it is
$\alpha\mapsto t\Ch_m(\alpha)$, where $t=\Td_{\leq m}(T_f)$.
We omit coefficients in $\Lambda$ from the vertices; the dashed arrows
indicate the three hidden edges.
\begin{equation}\label{eq:arr-cube}
\begin{tikzpicture}[
  x=1cm,y=.9cm,>=Stealth,
  every node/.style={font=\small},
  vertex/.style={inner sep=3pt,fill=white},
  edge/.style={->,thin},
  lab/.style={fill=white,inner sep=1.5pt},
  bridge/.style={preaction={draw=white,line width=3pt}}]
 \node[vertex] (X0) at (0,3.4) {$K_{\leq m}(X)$};
 \node[vertex] (Y0) at (6,3.4) {$K_{\leq n}(Y)$};
 \node[vertex] (X1) at (2,5.2) {$K_{\leq m}(X)$};
 \node[vertex] (Y1) at (8,5.2) {$K_{\leq n}(Y)$};
 \node[vertex] (C0) at (0,0) {$\CHle m(X)$};
 \node[vertex] (D0) at (6,0) {$\CHle n(Y)$};
 \node[vertex] (C1) at (2,1.8) {$\CHle m(X)$};
 \node[vertex] (D1) at (8,1.8) {$\CHle n(Y)$};
 \draw[edge,dashed] (C1) -- node[lab,left,pos=.64] {$\rho_X$} (X1);
 \draw[edge,dashed] (C1) -- node[lab,above,pos=.5] {$P_f$} (D1);
 \draw[edge,dashed] (C0) -- node[lab,below right,pos=.5] {$a^{-d}\Phi_X^a$} (C1);
 \draw[edge] (D1) -- node[lab,right] {$\rho_Y$} (Y1);
 \draw[edge] (C0) -- node[lab,left] {$\rho_X$} (X0);
 \draw[edge,bridge] (D0) -- node[lab,right,pos=.45] {$\rho_Y$} (Y0);
 \draw[edge,bridge] (X0) -- node[lab,above,pos=.55] {$f_!$} (Y0);
 \draw[edge] (X1) -- node[lab,above] {$f_!$} (Y1);
 \draw[edge] (X0) -- node[lab,above left] {$A_f^a$} (X1);
 \draw[edge] (Y0) -- node[lab,above left] {$\psi_Y^a$} (Y1);
 \draw[edge] (C0) -- node[lab,below] {$P_f$} (D0);
 \draw[edge] (D0) -- node[lab,below right] {$\Phi_Y^a$} (D1);
\end{tikzpicture}
\end{equation}
The upper face commutes by \eqref{eq:ARR}. Since the four
$\rho$-edges are isomorphisms, the lower face commutes as well:
$\Phi_Y^aP_f=P_f\,(a^{-d}\Phi_X^a)$, which is \eqref{eq:common-Adams}.
\end{proof}

\begin{lemma}\label{lem:Adams-Todd}
For every positive integer $a$ invertible
in $\Lambda$ one has, on $\CHle m(X)_{\Lambda}$,
\begin{equation}\label{eq:rho-Adams}
 \rho_X^{-1}\bigl(\theta^a(T_f^\vee)^{-1}\psi_X^a\bigr)\rho_X
       =a^{-d}\Phi_X^a.
\end{equation}
\end{lemma}

\begin{proof}
For a line bundle $L$, put $z=c_1^{\CH}(L)$. Then
$\Ch_m(\theta^a(L^\vee))=[\sum_{j=0}^{a-1}e^{-jz}]_{\leq m}$.
The universal series identity
\begin{equation}\label{eq:Todd-Adams-series}
 \frac{z}{1-e^{-z}}
 \left(\sum_{j=0}^{a-1}e^{-jz}\right)^{-1}
 =\frac{z}{1-e^{-az}}
 =a^{-1}\frac{az}{1-e^{-az}}
\end{equation}
follows from $\sum_{j=0}^{a-1}e^{-jz}=(1-e^{-az})/(1-e^{-z})$, and
has coefficients in $R_m[1/a]$ through degree $m$. Its left-hand side
is the series of the multiplicative class
$E\mapsto\Td_{\leq m}(E)\,\Ch_m(\theta^a(E^\vee))^{-1}$, while the
right-hand side is $a^{-1}$ times the Todd series evaluated at $az$;
on a line bundle, the latter corresponds to $a^{-1}\Phi_X^a$ applied
to the Todd class. Taking multiplicative extensions to $T_f$, a virtual
bundle of rank $d$, we obtain
\begin{equation}\label{eq:Chow-Adams-Todd}
 t\,\Ch_m\bigl(\theta^a(T_f^\vee)^{-1}\bigr)=a^{-d}\Phi_X^a(t).
\end{equation}
This is the truncated Todd--Adams identity; for its rational
predecessor, see \cite[Ch.~III, Lemma~4.2]{FultonLang}.
For $u\in\CHle m(X)_{\Lambda}$, multiplicativity of $\Ch_m$ and
\eqref{eq:ch-adams} now give
\[\begin{aligned}
 \rho_X^{-1}\bigl(\theta^a(T_f^\vee)^{-1}\psi_X^a\bigr)\rho_X(u)
 &=t\,\Ch_m\bigl(\theta^a(T_f^\vee)^{-1}\bigr)\Phi_X^a(t^{-1}u)\\
 &=a^{-d}\Phi_X^a(t)\Phi_X^a(t^{-1}u)\\
 &=a^{-d}\Phi_X^a(u).
\end{aligned}\]
The last equality uses that $\Phi_X^a$ is a ring map, so
$\Phi_X^a(t^{-1})=\Phi_X^a(t)^{-1}$.
\end{proof}

\section*{Funding}
Both authors were supported by FEDER/AEI/MICINN under Grant No.~PID2022-142024NB-I00 (SIGA).

\section*{Use of artificial intelligence}
We used ChatGPT (OpenAI) to discuss the proof and to draft and revise
the exposition. We also used Fable~5.2 to review a late draft.
We take full responsibility for the mathematical arguments and the final text.

\begingroup
\footnotesize
\setstretch{1}

\endgroup
\par\bigskip
\begin{flushleft}
\footnotesize
\setstretch{1}
\paperemails
\end{flushleft}
\end{document}